\documentclass[12pt]{article}

\usepackage[utf8]{inputenc}
\usepackage{graphicx}
\usepackage{xcolor}
\usepackage{natbib}
\usepackage[toc,page]{appendix}
\usepackage{geometry}
\usepackage{microtype}
\usepackage{setspace}
\usepackage{fullpage}
\usepackage{lineno}
\usepackage{amsmath, amssymb, amsthm}
\usepackage{graphics}
\usepackage{xcolor}

\usepackage{amsmath, amssymb, amsthm}
\theoremstyle{plain}
\newtheorem{lemma}{Lemma}
\newtheorem{theorem}{Theorem}
\newtheorem{proposition}{Proposition}
\newtheorem{corollary}{Corollary}
\theoremstyle{definition}

\newtheorem{example}{Example}
\theoremstyle{remark}

\newcommand{\E}{\mathsf{E}}

\newcommand{\var}{\mathrm{var}}
\newcommand{\cov}{\mathrm{cov}}

\newcommand{\pr}{\mathrm{P}}

\DeclareMathOperator{\vech}{vech}

\DeclareMathOperator{\bdiag}{bdiag}
\newcommand{\tr}{\operatorname{tr}}
\newcommand{\tsp}{\mathrm{\scriptscriptstyle T}}
\newcommand{\rN}{\mathrm{N}}
\newcommand{\emax}{\gamma_{\mathrm{max}}}
\newcommand{\emin}{\gamma_{\mathrm{min}}}
\newcommand{\R}[1]{\mathbb{R}^{#1}}

\title{Uniform inference in linear mixed models}
\author{Karl Oskar Ekvall$^{\star, \dagger}$ \quad Matteo Bottai$^{\dagger}$ \\
{\normalsize $^\star$Department of Statistics, University of Florida} \\
{\normalsize $^\dagger$Division of Biostatistics, Institute of Environmental
Medicine, Karolinska Institutet}\\
{\tt k.ekvall@ufl.edu \quad
matteo.bottai@ki.se}}
\date{\normalsize \today}
\date{\today} 

\usepackage{hyperref}
\usepackage{lipsum}

\begin{document}

\onehalfspacing

\maketitle

\begin{abstract}
    We provide finite-sample distribution approximations, that are uniform in the parameter, for inference
    in linear mixed models. Focus is on variances and
    covariances of random effects in cases where existing theory fails because
    their covariance matrix is nearly or exactly singular, and hence near or at
    the boundary of the parameter set. Quantitative bounds on the differences
    between the standard normal density and those of linear combinations of the
    score function enable, for example, the assessment of sufficient sample
    size. The bounds also lead to useful asymptotic theory in settings where
    both the number of parameters and the number of random effects grow with the
    sample size. We consider models with independent clusters and ones with
    a possibly diverging number of crossed random effects, which are notoriously
    complicated. Simulations indicate the theory leads to practically relevant
    methods. In particular, the studied confidence regions, which are
    straightforward to implement, have near-nominal coverage in finite samples
    even when some random effects have variances near or equal to zero, or
    correlations near or equal to $\pm 1$.
\end{abstract}

\section{Introduction} \label{sec:intro}

Linear mixed models, and random effects in particular, are used routinely to
model dependence and effect heterogeneity. However, while random effects are
convenient for specifying a model, they often complicate inference. For example,
it is well known that, when testing if the variance of a random effect is zero,
common test-statistics have nonstandard asymptotic distributions because the
parameter is on the boundary of the parameter set; see for example
\citet{self1987asymptotic} and \citet{geyer1994asymptotics} and the references therein. In
general, the distributions, and hence the appropriate critical values for tests,
depend on the particular boundary point, the structure of the parameter set, and
the test-statistic. By contrast, at any fixed interior point of the
parameter set, test-statistics such as score, Wald, and likelihood ratio all
have asymptotic chi-squared distributions under regularity conditions. It is common, therefore, to use chi-squared quantiles as critical values for
interior points. Unfortunately, doing so does not lead to reliable
inference in general. Confidence regions obtained by inverting the tests often
have uniform coverage probabilities quite different from nominal, even
asymptotically. The regions are often overly conservative, but they can also be
invalid. To date, the issues have been addressed only in a few special cases of
the settings considered here. Importantly, existing results on near-boundary
inference preclude many dependence structures common in mixed models.

Let $X \in \R{n \times p}$ and $Z \in \R{n\times
q}$ be non-stochastic matrices corresponding to a vector $\beta \in \R{p}$ of
fixed effects and a vector $U \in \R{q}$ of random effects, respectively.
Suppose $U$ is multivariate normal with mean zero and covariance matrix $\Psi
\in \R{q\times q}$, and that a vector $Y\in \R{n}$ of responses satisfies
\begin{equation} \label{eq:model} 
  Y = X \beta + Z U + E, 
\end{equation}
where $E\sim \rN(0, \psi_r I_n)$ and $U \sim \rN(0, \Psi)$, independently, and
$\Psi$ is parameterized by the first $r - 1$ elements of the vector $\psi =
[\psi_1, \dots, \psi_r]^\tsp$; see Section \ref{sec:model} for details. A main
goal is reliable inference on $\psi$, and in particular confidence regions with
good coverage properties. Focus is on parameters near or at the
boundary of the parameter set, because those are often of interest in practice,
and many common methods fail near the boundary. The boundary includes points
where $\Psi$ is singular, which happens, for example, if some random effects
have perfect correlation or vanishing variances. Thus, to reliably assess the practical significance of random effects, confidence regions need to have good coverage properties near the boundary.

\begin{figure}
  \includegraphics[width=\linewidth]{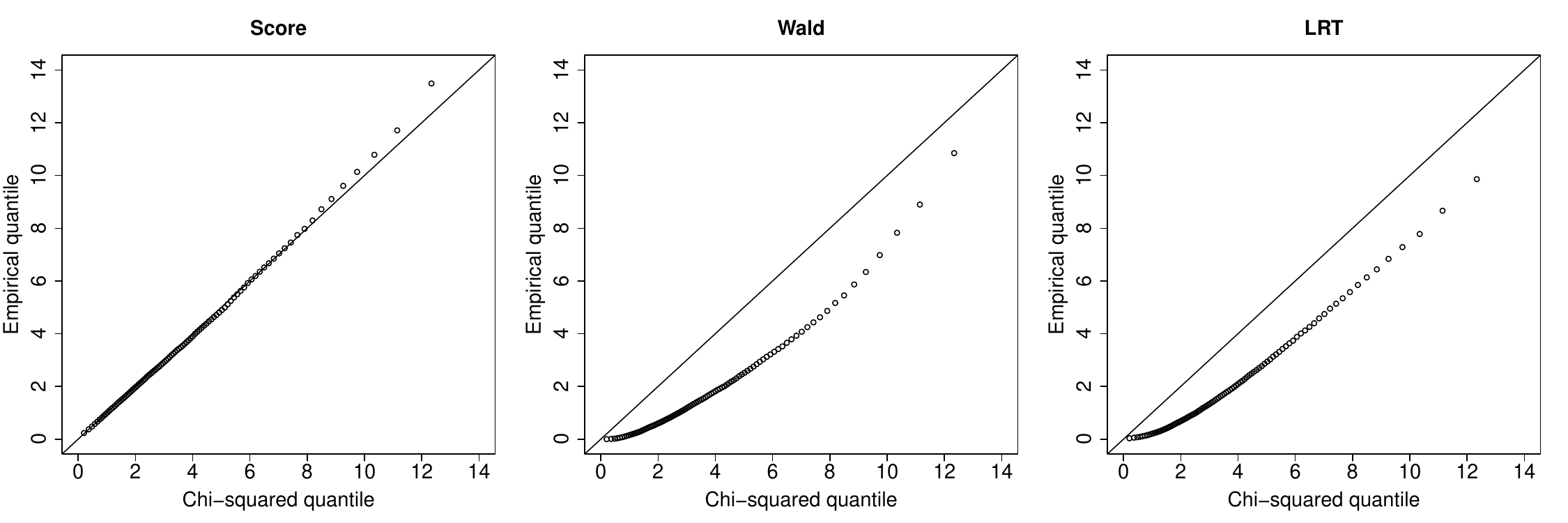}
\caption{Quantiles of test-statistics evaluated at a true $\psi \in \R{4} $ that is near the
boundary. For $i \in \{1, \dots, 50\}$, $Y_{i} = Z_i U_i + E_i$, where $U_i\sim
\rN(0, \Psi_1)$, $\Psi_1 \in \R{2\times 2}$ has diagonal elements $\psi_1 =  \psi_3 = 10^{-3}$ and off-diagonal $\psi_2 = 0$, and $E_i \sim \rN(0, \psi_4)$, $\psi_4 = 1$. Elements of
$Z_i \in \R{5\times 2}$ were drawn prior to simulations from a Bernoulli distribution with mean $1/2$. Empirical quantiles based on 10 000
replications.}
\label{fig:intro}
\end{figure}

A key issue is that the distributions of test-statistics in general depend not
only on whether a parameter is on the boundary or not, but on its proximity to
the boundary
\citep{moran1971maximumlikelihood,rotnitzky2000likelihoodbased,stern2000likelihood,bottai2003confidence,crainiceanu2004likelihood,ekvall2022confidence,zhang2025fast}.
Figure \ref{fig:intro} illustrates this in a special case of \eqref{eq:model}.
Clearly, even though the parameter is interior, the distributions of Wald and
likelihood ratio statistics are quite different from chi-squared. Moreover, the
distributions are different from each other, as also noted in other recent
research \citep{battey2024anomaly}. By contrast, the figure shows the score
statistic is approximately chi-squared distributed. Evidently, it is
inappropriate to use the same critical values for all test-statistics, even at
interior points. These nuances are not reflected in pointwise asymptotic theory,
but they are revealed by uniform results (Section \ref{sec:boundary}). 

Suppose for now that $\beta = 0$ is known and let $S(\psi) \in \R{r}$ be the
score, or the gradient of the log-likelihood at $\psi$, and $\mathcal{I}(\psi) =
\cov_\psi\{S(\psi)\} \in \R{r\times r}$ the Fisher information matrix. Here
and elsewhere, the subscript $\psi$ indicates the covariance is with respect to
the distribution indexed by $\psi$. Let also $W^S(\psi) =
\mathcal{I}(\psi)^{-1/2}S(\psi)$ and, for non-zero $v \in \R{r}$, let $g(\cdot;
v, \psi)$ be the density of $v^\tsp W^S(\psi)$ under $\psi$; that is, when
$\psi$ is the true parameter. Our main finite-sample results in Section
\ref{sec:asy_dist} upper bound $\vert g(t; v, \psi) - \phi(t)\vert$, where
$\phi(\cdot)$ is the standard normal density. The bounds are uniform over
suitably chosen sets of $t$, $v$, and $\psi$, including ones with boundary
points. Upon letting $n$ tend to infinity, possibly along with $r$ and $p$, the
bounds also lead to uniform asymptotic results.

We treat both settings with independent clusters and ones with crossed random
effects. The latter are particularly challenging and not handled by existing
theory for inference near the boundary. Moreover, to the best of our knowledge,
there are no existing finite-sample results similar to ours. In fact, the first
results on pointwise asymptotic normality of maximum likelihood estimators of
interior parameters of fixed dimension only appeared shortly before submission
of our work \citep{jiang2024precise,lyu2024increasing}. However, while impressive and
useful for other purposes, those results are not uniform, and cannot be made so
near the boundary. Consequently, they do not lead to reliable confidence regions
for $\psi$ in general; see Section \ref{sec:boundary} for details and an example illustrating the limitations of pointwise convergence in distribution.

To be sure, we are not the first to consider score-based inference in mixed
models \citep[see for
example][]{verbeke2003use,qu2013linear,zhu2006generalized}. However, we are
aware of no results similar to ours. \citet{zhang2025fast} developed asymptotic
theory for score-based confidence regions for a variance component, but only for
the special case of \eqref{eq:model} where $\Psi = \psi_1 I_q$; that is, with a
single variance component. That setting is substantially simpler, partly because
there is only one variance component, but also because observations are
independent after transforming the data by the left singular vectors of $Z$.
Thus, theory can be developed assuming independence, which is not the case in
general under \eqref{eq:model}. Indeed, dependence is often complicated,
especially with crossed random effects, and there can be few independent vectors
even as $n$ grows \citep[see for example][]{sung2007monte,greven2008restricted}.
\citet{ekvall2022confidence} provided asymptotic results in a related setting, but
with a diverging number of independent observations, a fixed number of
parameters, diagonal $\Psi$, and singular Fisher information, neither of which
is true here in general.

\section{Background on inference near the boundary} \label{sec:boundary}

We will often first state results where $\beta = 0$ is known, so that $\psi$ is
the only parameter. This lets us focus on the key issues, which arise when
making inferences about the covariance parameters, and it simplifies notation.
Then we address settings where $\beta$ is unknown.

Let $T(\psi) = T(\psi; Y, X, Z)$ denote a generic test-statistic, defined for
$\psi$ in some parameter set $\mathbb{P}$. Define a $(1 - \alpha) \in (0,
1)$ confidence region obtained by inverting $T$ by
\begin{equation} \label{eq:region}
  \mathbb{C}(\alpha) = \{\psi \in \mathbb{P}: T(\psi) \leq q_{1 - \alpha}(\psi)\}.
\end{equation}
If the critical value $q_{1 - \alpha}(\psi)$ is the $(1-\alpha)$th quantile of
the exact distribution of $T(\psi)$, for every $\psi \in \mathbb{P}$, then the
confidence region has uniformly correct coverage probability. That is, $\sup_{\psi \in
\mathbb{P}} \vert\pr_\psi\{\psi \in \mathbb{C}(\alpha)\} - (1 - \alpha)\vert =
0$. However, because exact distributions are typically unavailable, we instead
look for useful approximations.Assume for simplicity, for the remainder of Section \ref{sec:boundary}, that
$\mathbb{P}$ does not depend on $n$. Then we would like to pick $q_{1 -
\alpha}(\psi)$ such that, for every compact $\mathbb{A} \subseteq \mathbb{P}$,
\begin{equation} \label{eq:asy_uni_cover}
  \lim_{n\to \infty} \sup_{\psi \in \mathbb{A}} \vert\pr_\psi\{\psi \in 
    \mathbb{C}_n(\alpha)\} - (1 - \alpha)\vert = 0.
\end{equation}
The added subscript $n$ indicates the confidence region depends on $n$; we omit these subscript when stating finite-sample results. When \eqref{eq:asy_uni_cover} holds, we say $\mathbb{C}_n(\alpha)$ has
asymptotically correct coverage probability on compact sets. We will see that
the compactness restriction can sometimes be relaxed.

The following lemma, which is not specific to mixed models, gives a useful
characterization of \eqref{eq:asy_uni_cover}. The lemma is somewhat similar to,
but applies more generally than, Lemma 2.5 in \citet{ekvall2022confidence}. Proofs are in the Supplementary Material unless otherwise stated.

\begin{lemma} \label{lem:equiv_conv}
  Equation \eqref{eq:asy_uni_cover} holds for every compact $\mathbb{A}
  \subseteq\mathbb{P}$ if and only if, for every sequence
  $(\psi_n)$ convergent in $\mathbb{P}$,
  \begin{align} \label{eq:seq_prob_cov}
    \lim_{n\to \infty}\vert\pr_{\psi_n}\{\psi_n \in 
  \mathbb{C}_n(\alpha)\} - (1 - \alpha)\vert = 0.
  \end{align}
\end{lemma}

We will use that \eqref{eq:seq_prob_cov} implies \eqref{eq:asy_uni_cover} in our
proofs, but the other direction of the equivalence is also important: It tells
us that pointwise convergence in distribution cannot be used to establish
\eqref{eq:asy_uni_cover} in general--one must consider sequences of parameters.
The next result, which essentially says \eqref{eq:seq_prob_cov} holds if
a test-statistic has the same asymptotic distribution under any convergent sequence
of parameters, is also not specific to mixed models; see for example
\citet{mikusheva2007uniform} for a similar result.

\begin{lemma} \label{lem:same_asy_dist}
  Suppose that, for a test-statistic $T_n$ and continuous cumulative distribution function $F$, it holds
  for every sequence $(\psi_n)$ convergent in $\mathbb{P}$ and $t \in \R{}$,
  that $F_n(t) = \pr_{\psi_n}\{T_n(\psi_n) \leq t\} \to F(t)$ as $n\to \infty$.
  Then, for any $\alpha \in (0, 1)$, $\mathbb{C}_n(\alpha)$ defined by
  \eqref{eq:region} with $T = T_n$ and $q_{1 - \alpha}(\psi) = F^{-}(1 - \alpha) = \min\{t:
  F(1 - \alpha) \geq t\}$, satisfies \eqref{eq:seq_prob_cov}.
\end{lemma}

In later sections we will verify the conditions of Lemma \ref{lem:same_asy_dist}
for score-based confidence regions. Conversely, the
conditions of the lemma do not hold for likelihood ratio and Wald statistics in
general. To see this it suffices to consider constant sequences
with $\psi_n = \psi$, which are trivially convergent. Under such sequences,
classical results say likelihood ratio and Wald statistics have different
asymptotic distributions depending on whether $\psi$ is an interior or boundary
point. This does not imply that those statistics cannot give a confidence
region satisfying \eqref{eq:asy_uni_cover}, but Lemma \ref{lem:same_asy_dist}
does not apply to them, so $q_{1 - \alpha}(\psi)$ would have to
depend on $\psi$ in some nontrivial way in general. It is not enough to
use different critical values for boundary and interior points, as the following
example illustrates. We have made the example as simple as possible
while still illustrating key issues to be addressed more generally.

\begin{example} \label{ex:intro}
Suppose that, independently for $i \in \{1, \dots, n\}$,
\begin{equation} \label{eq:ex_model}
  Y_i \sim \rN(0, 1 + \psi_1).
\end{equation}
This is a special case of \eqref{eq:model} with known $\beta = 0$, $Z = I_n$,
$\Psi = \psi_1 I_n$, and $\psi \in \mathbb{P} = [0, \infty) \times \{1\}$.
Because $\psi_2 = 1$, we simplify notation and write $\psi = \psi_1 \in [0,
\infty)$ for the remainder of the example. Up to a constant, the log-likelihood is $\ell_n(\psi) =
-n\{\log(1 + \psi) + M_n/(1 + \psi)\}/2$, where $M_n = \sum_{i = 1}^n Y_i^2/n$.
The score is $S_n(\psi) = \nabla \ell_n(\psi) = -n\{(1 + \psi)^{-1} - M_n(1 +
\psi)^{-2}\}/2$ and the Fisher information is $\mathcal{I}_n(\psi) =
\var_{\psi}\{S_n(\psi)\} = n / \{2(1 + \psi)^2\}$. Define score and Wald
statistics by, respectively,
\[
  W_n^{S}(\psi) = \mathcal{I}_n(\psi)^{-1/2} S_n(\psi) =
    (n/2)^{1/2}\{M_n/(1 + \psi) - 1\}
\]
and $W_n^W(\psi) = (\hat{\psi}_n - \psi) \mathcal{I}_n(\psi)^{1/2}$, where
$\hat{\psi}_n = \max(M_n - 1, 0)$ is the maximum likelihood estimator. The
likelihood ratio test-statistic is $T_n^{L}(\psi) = 2\{\ell_n(\hat{\psi}_n) -
\ell_n(\psi)\}$. It is clear that $\hat{\psi}_n$ can be zero with substantial
probability. More specifically, since $n M_n/(1 + \psi) \sim \chi^2_n$, a normal
approximation gives $\pr_\psi(\hat{\psi} = 0)
\approx \Phi\{-(n/2)^{1/2}\psi /(1 + \psi)\}$, where $\Phi$ is the
standard normal cumulative distribution function. This probability is $1/2$ if
$\psi = 0$, and is close to that if $n^{1/2}\psi$ is small. This
suggests a normality approximation for the maximum likelihood estimator may
perform poorly near the boundary. Intuitively, then, a chi-squared approximation
for the likelihood ratio statistic may also be inappropriate. The following
result formalizes this intuition.

\begin{proposition} \label{prop:example}
Suppose \eqref{eq:ex_model} holds with
$\psi = \psi_n \in [0, \infty)$ satisfying $\psi_n n^{1/2} \to a \in [0,
\infty]$ as $n\to \infty$. Then, with $W_1 \sim \rN(0, 1)$, in distribution,
\begin{align}
  W_n^S(\psi_n)& \to W_1, \label{eq:score_conv_ex} \\ 
  W_n^W(\psi_n) &\to \max(W_1, -a 2^{-1/2}),\label{eq:wald_conv_ex} \\
  T_n^L(\psi_n)& \to 2W_1 \max(W_1, -a 2^{-1/2}) - \max(W_1, -a2^{-1/2})^2.
  \label{eq:lrt_conv_ex}
\end{align}
\end{proposition}
Results similar to Proposition \ref{prop:example} are given by \citet{moran1971maximumlikelihood}
and \citet{stern2000likelihood}. Clearly, the test-statistics are not asymptotically equivalent under sequences of parameters. Statistics which use the maximum
likelihood estimator behave irregularly since there is an appreciable
probability the estimator is on the boundary, as is illustrated by
\eqref{eq:wald_conv_ex}--\eqref{eq:lrt_conv_ex}.

Proposition \ref{prop:example} implies a confidence region based on
$T_n^S(\cdot) = \{W_n^S(\cdot)\}^2$ satisfies \eqref{eq:asy_uni_cover}. Indeed, by
\eqref{eq:score_conv_ex}, Lemma \ref{lem:same_asy_dist} applies with $T_n = T_n^S$ upon
taking $F$ to be the cumulative distribution function of $\chi^2_1$. We will see
that this generalizes to more complex settings.

Suppose $a = 0$ in Proposition \ref{prop:example}, including $\psi_n = 0$ for
every $n$ as a special case. Then we recover familiar asymptotic distributions
for inference on boundary points. The asymptotic distribution of
$T_n^L(\psi_n)$ simplifies to that of $\max(W_1, 0)^2$, a mixture of $\chi^2_1$
and a point mass at zero. Conversely, if $a = \infty$, which happens for
example if $\psi_n = \psi > 0$, we recover the classical result that
$T_n^S(\psi)$, $T_n^W(\psi) = \{W_n^W(\psi)\}^2$, and $T_n^L(\psi)$ are all
asymptotically $\chi^2_1$.

To see why a likelihood ratio confidence region can be unreliable
in our settings, let $\mathbb{C}_n^L(\alpha)$ be defined by \eqref{eq:region} with $T = T_n^L$ and, for interior $\psi > 0$, $q_{1 - \alpha}(\psi) = c_{1, 1 -
\alpha}$, the $(1- \alpha)$th quantile of $\chi^2_1$. Then, regardless of which critical value is used for the boundary point $\psi = 0$, we get by \eqref{eq:lrt_conv_ex} with $\psi_n = 1/n$,
\begin{align*}
  \sup_{\psi \in [0, 1]}
  \vert\pr_{\psi}\{\psi \in \mathbb{C}^L_n(\alpha)\} - (1 - \alpha)\vert &\geq
  \vert \pr_{1/n}\{T_n^L(1/n) \leq c_{1, 1 - \alpha}\} - (1 - \alpha)\vert \\
  & \to \vert \pr\{\max(W_1, 0)^2 \leq c_{1, 1 - \alpha}\} - (1 - \alpha)\vert,
\end{align*}
which is greater than zero. For example, with $\alpha = 0.05$, $\pr\{\max(W_1,
0)^2 \leq c_{1, 0.95}\} = \pr(W_1 \leq 1.96) = 0.975 \neq 0.95$. Thus,
$\mathbb{C}_n^L(\alpha)$ does not satisfy \eqref{eq:asy_uni_cover}.
Intuitively, this happens because $T_n^L(\psi)$ behaves as if $\psi$ is on the
boundary when it is close enough to it. Similar arguments apply to the Wald
statistic. For reasons illustrated by this example, we shall focus on
score-statistics.
  
\end{example}

\section{Model and test-statistic}\label{sec:model}

Equation \eqref{eq:model} implies
\begin{align} \label{eq:marginal}
  Y \sim \rN(X\beta, \Sigma), ~~\Sigma = Z \Psi Z^\tsp +  \psi_r I_n,
\end{align}
where $\Psi \in \R{q\times q}$ is parameterized by $\psi_{-r} = [\psi_1, \dots,
\psi_{r-1}]^\tsp$. We assume each element of $\Psi$ is either known to be zero
or equal to one of the elements of $\psi_{-r}$. Thus, the elements of
$\psi_{-r}$ are variances or covariances of random effects. To simplify
notation, we consider $\Psi = \Psi(\psi)$ a function of $\psi$ and let $H_j = \partial \Psi / \partial \psi_j$, $j \in \{1,
\dots, r\}$.  The elements of $H_j$ are all zeros and
ones, with at least one non-zero element for each $j < r$; $H_r = 0$; and
$\Psi(\psi) = \sum_{j = 1}^r \psi_{j} H_j$. The parameter set for $\psi$ is
$\mathbb{P} = \{\psi \in \R{r}: \Psi(\psi) \geq 0, \psi_r > 0\}$, where
$\Psi(\psi) \geq 0$ indicates positive semi-definiteness. For example, in models
with a single random effect, $\Psi(\psi) = \psi_1 I_q$, $H_1 = I_q$, and
$\mathbb{P} = [0, \infty) \times (0, \infty)$. We often omit the
argument to $\Psi$ for simplicity. Similarly, we write $\Sigma(\psi) = Z
\Psi(\psi)Z^\tsp + \psi_r I_n$ when the argument needs to be emphasized, and omit it otherwise.

The log-likelihood for $\theta = (\beta, \psi) \in \R{p}\times \mathbb{P}$
corresponding to \eqref{eq:marginal} is, up to a constant,
\[
  \ell(\theta) = \ell(\theta; Y, X, Z) = -\frac{1}{2}\{\log \vert \Sigma(\psi)\vert +
  (Y - X\beta)^\tsp \Sigma(\psi)^{-1}(Y - X\beta)\}.
\]
We abuse notation somewhat and define score functions for $\theta$, $\beta$, and $\psi$, respectively, by
\begin{align*}
  S(\theta) = \frac{\partial \ell(\theta)}{\partial \theta} \in \R{p + r},~~
  S(\beta; \psi) = \frac{\partial \ell(\beta, \psi)}{\partial \beta} \in
\R{p},~~
S(\psi; \beta) = \frac{\partial \ell(\beta, \psi)}{\partial \psi} \in
\R{r}.
\end{align*}
Similarly, $S(\psi_j) = \partial \ell_n(\beta,
\psi)/\partial \psi_j$, $j  \in \{1, \dots, r\}$. Let $A_{j} = \Sigma^{-1/2}(\partial \Sigma / \partial \psi_j)\Sigma^{-1/2}$, so $A_{j} = \Sigma^{-1/2}ZH_jZ^\tsp \Sigma^{-1/2}$ for $j \in \{1,
\dots, r - 1\}$, and $A_r = \Sigma^{-1}$. Then, with
$R = \Sigma^{-1/2}(Y - X\beta)$,
\begin{align} \label{eq:score_psij}
  S(\psi_j) &= \frac{1}{2}\{R^\tsp A_{j}
    R - \tr(A_{j})\},~~ j \in \{1, \dots, r\},
\end{align}
and $S(\beta) = X^\tsp \Sigma^{-1/2}R$. The Fisher information matrix is $\mathcal{I}(\theta) = \cov_{\theta}\{S(\theta)\}$. Since the scores for the
$\psi_j$ are quadratic forms in $R$ and the score for $\beta$ is linear in $R$,
an expression for $\mathcal{I}(\theta)$ follows from the following routine
result. The result does not require normality.
\begin{lemma} \label{lem:finf}
   Suppose $A_1, A_2 \in \R{n\times n}$ are symmetric and that $R \in \R{n}$ is a
   random vector with mean zero and identity covariance matrix. Then $\E(R^\tsp
   A_1 R) = \tr(A_1)$. If in addition $\E(R_i^3) = 0$ for all $i$, then
   $\E\{(R^\tsp A_1 R) R\} = 0$; and if also $\E(R_i^4) = 3$ for all $i$,
   then $\cov(R^\tsp A_1 R, R^\tsp A_2 R) = 2\tr(A_1 A_2)$.
  \end{lemma}
Lemma \ref{lem:finf} is essentially well known, but a proof is in
the Supplementary Material for completeness. The expressions in the proof give $\cov_{\theta}\{S(\theta)\}$ also in settings where the model is
misspecified so that some conditions in Lemma \ref{lem:finf} do not hold, for
example when $Y$ is not normally distributed. It follows from the lemma that
$\mathcal{I}(\theta)$ is block-diagonal with leading $p\times p$ block and
trailing $r\times r$ block given by, respectively,
\begin{equation} \label{eq:finf}
  \mathcal{I}(\beta; \psi) = X^\tsp \Sigma^{-1}X,~~ \mathcal{I}_{ij}(\psi) 
    = \tr(A_iA_j)/2, ~~ i,j \in \{1, \dots, r\}.
\end{equation}
As the notation suggests, the information for $\psi$ does not depend on $\beta$.
To state the next result, let $\Psi(v) = \sum_{j = 1}^r v_j H_j$ and $\Sigma(v)
= Z \Psi(v)Z^\tsp + v_r I_n$ for any $v \in \R{r}$. Both $\Psi(v)$ and
$\Sigma(v)$ are symmetric, but they need not be positive semi-definite. Let also
$\mathbb{S}^{r-1} = \{v\in \R{r}:\Vert v\Vert = 1\}$.

\begin{theorem} \label{thm:finf_pd}
  The matrix $\mathcal{I}(\theta)$ is positive definite if and only if
  $\mathcal{I}(\beta; \psi)$ and $\mathcal{I}(\psi)$ are. The matrix
  $\mathcal{I}(\beta; \psi)$ is positive definite if and only if $X$ has full
  column rank $p \leq n$, and $\mathcal{I}(\psi)$ is positive definite if and
  only if, for every $v\in \mathbb{S}^{r - 1}$,
  \begin{equation} \label{eq:pdinf_cond}
   \Vert \Sigma(v)\Vert > 0.
  \end{equation}
\end{theorem}
The condition involving \eqref{eq:pdinf_cond} can only hold if $r \leq n(n + 1)/2$, and $r$ is typically substantially smaller than that. It is almost trivial to show that the conditions in
Theorem \ref{thm:finf_pd} ensure identifiability. However, identifiability does not imply a positive definite information matrix in general \citep{rothenberg1971identification}.
Indeed, if \eqref{eq:model} were parameterized in terms of a matrix $\Lambda \in
\R{r\times r}$ such that $\Psi = \Lambda \Lambda^\tsp$, which can be done
without losing identifiability, then the information matrix would be
singular at points
where $\Lambda$ is singular \citep{ekvall2022confidence, guedon2024bootstrap}; see also
\citet{chesher1984testing,cox2000theoretical,lee1986specification}. By contrast, since the
conditions in Theorem \ref{thm:finf_pd} do not say anything about $\theta$, here $\mathcal{I}(\theta)$ is either positive definite
for every $\theta$ or for none. Consequently, here, singular information implies an unidentifiable parameter.

The following corollary is useful in some examples.

\begin{corollary} \label{corol:psi_in_pd}
  If $Z$ has full column rank $0 < q < n$, then $\mathcal{I}(\psi)$ is
  positive definite.
\end{corollary}

The condition in Corollary \ref{corol:psi_in_pd} is not necessary.
For example, the proof reveals that if $\Psi = \psi_1
I_q$, $\mathcal{I}(\psi)$ is positive definite unless $ZZ^\tsp$ is proportional
to the identity. Settings with $\Psi(\psi) = \psi_1 I_q$ are fairly common in
applications, and were studied by \citet{zhang2025fast}. More
generally, a diagonal $\Psi$ leads to a variance components model where 
$\Sigma = \sum_{j = 1}^r \psi_j K_j$, $K_r = I_n$, and $K_j$ positive semi-definite for $j < r$. For that model, Theorem \ref{thm:finf_pd} says $\mathcal{I}(\psi)$ is positive definite if and only if $K_1, \dots, K_r$ are linearly independent.

When $\mathcal{I}(\theta)$ is invertible, define
\begin{equation} \label{eq:score_stat}
  W^S(\theta) = \mathcal{I}(\theta)^{-1/2}S(\theta),~~ 
  T^S(\theta) = \Vert W^S(\theta)\Vert^2,
\end{equation}
where $\Vert \cdot \Vert$ is the Euclidean norm. Unlike Wald statistics, for
example, $T^S$ is invariant under differentiable reparameterizations with full
rank Jacobian. Consequently, our results are not specific to the considered
parameterization.

For inference on $\psi$ only it is common to use the restricted likelihood
\citep{patterson1971recovery}.  It often gives estimators that are less biased,
and less likely to be on the boundary \citep{stern2000likelihood}. Suppose $X$
has column rank $p < n$, and let $V \in\R{n\times (n-p)}$ be a semi-orthogonal
matrix such that $V^\tsp X = 0$. Let also $\tilde{Y} = V^\tsp Y$, $\tilde{Z} =
V^\tsp Z$, and $\tilde{\Sigma} = V^\tsp \Sigma V$. Then
\begin{equation} \label{eq:model_reml}
  \tilde{Y} \sim \rN(0, \tilde{\Sigma}), ~~ \tilde{\Sigma}= \tilde{Z} \Psi \tilde{Z}^\tsp + \psi_r I_{n - p}.
\end{equation}
The restricted likelihood is the likelihood for \eqref{eq:model_reml}. Having
not specified a particular $Z$ or $X$, \eqref{eq:model_reml} is essentially a
special case of \eqref{eq:marginal} with known $\beta = 0$ and sample size $n -
p$. Therefore, many results for the usual likelihood apply to the restricted
likelihood after minor changes to notation. For this reason, we often first
state results assuming \eqref{eq:marginal} with known $\beta = 0$ and then apply
these to the restricted likelihood.

We write $\tilde{S}(\psi)$ and
$\tilde{\mathcal{I}}(\psi)$ for the restricted score and information matrix,
respectively. Using the similarity of 
\eqref{eq:marginal} and \eqref{eq:model_reml}, explicit expressions are immediate from \eqref{eq:score_psij}--\eqref{eq:finf}. In particular, $R$ and $A_j$, are replaced by, respectively, $\tilde{R} = \tilde{\Sigma}^{-1/2}\tilde{Y}$ and $\tilde{A}_j =
\tilde{\Sigma}^{-1/2} \tilde{Z} H_j \tilde{Z}^\tsp \tilde{\Sigma}^{-1/2}$, if $j
< r$, or $\tilde{A}_r = \tilde{\Sigma}^{-1}$. Similarly, $\tilde{W}^S(\psi) =
\tilde{\mathcal{I}}(\psi)^{-1/2}\tilde{S}(\psi)$ and $\tilde{T}^S(\psi) = \Vert
\tilde{W}^S(\psi)\Vert^2$. The
Supplementary Material contains additional remarks how inference based on the
different likelihoods can be implemented in practice.

\section{Approximate distributions} \label{sec:asy_dist}
\subsection{Finite-sample bounds and asymptotic normality}
We now turn to finite-sample distribution approximations for score statistics. Recall $A_j(\psi) = \Sigma(\psi)^{-1/2}\{\partial \Sigma(\psi) / \partial \psi_j\}\Sigma(\psi)^{-1/2}$, $j
\in \{1, \dots, r\}$. For $v \in \R{r}$ and $\psi \in \mathbb{P}$, let 
\[
  A(v, \psi) = \sum_{j = 1}^r v_j A_j(\psi) = \Sigma(\psi)^{-1/2}\Sigma(v)\Sigma(\psi)^{-1/2}.
\]
For $v$ and $\psi$ such that $\Vert A(v, \psi) \Vert_F > 0$, define also
\[
  a(v, \psi) = \frac{\Vert A(v, \psi)\Vert}{\Vert A(v, \psi)\Vert}_F \leq 1,
\]
where $\Vert \cdot \Vert$ and $\Vert \cdot \Vert_F$ are the spectral and
Frobenius norms, respectively. Recall, for symmetric matrices, the squared Frobenius norm is the sum of squared eigenvalues, and the spectral norm is the largest absolute eigenvalue. For $v$ and $\psi$ such that $\Vert A(v, \psi) \Vert_F = 0$, define $a(v, \psi) = 1$.  Note $\Vert A(v, \psi)\Vert_F > 0$ if and only if
$\Vert \Sigma(v)\Vert > 0$ since $\emin\{\Sigma(\psi)\}\geq \psi_r > 0$, where
$\emin(\cdot)$ is the smallest eigenvalue. Thus, by Theorem \ref{thm:finf_pd}
and the discussion following it, $\Vert A(v, \psi) \Vert_F > 0$ for all non-zero
$v \in \R{r}$ if $\mathcal{I}(\psi)$ positive definite and $\psi$ identifiable. Shortly, we will see the score is close to normal if $a(v, \psi)$ is small for appropriate $v$ and $\psi$. For intuition, recall that in Example \ref{ex:intro}, $\Sigma(v) = (1 + v) I_n$, and hence $a(v, \psi) = n^{-1/2}$. In more complicated examples given later there is an interplay between the identifiability of $\psi$, controlled by $\Sigma(v)$, and the regularity of $\Sigma(\psi)$.

The following result builds on a normal approximation for quadratic forms due to
\citet{zhang2025fast}. Recall $g(\cdot; v, \psi )$ is the density of $v^\tsp
W^S(\psi)$ under $\psi$ when \eqref{eq:marginal} holds with known $\beta = 0$,
so that $\theta = \psi$, and $\phi(\cdot)$ is the standard normal density.

\begin{lemma} \label{lem:zhang-ish}
  For any $t\in \R{}$,  $\psi \in \mathbb{P}$, and $v\in \mathbb{S}^{r - 1}$
  such that $\tilde{v} = \mathcal{I}^{-1/2}(\psi) v \neq 0$ and $a(\tilde{v},
  \psi)^2 < 1/8$, it holds that
  \begin{align} \label{eq:zhang-ish}
    \vert g(t; v, \psi) - \phi(t)\vert \leq  0.14\left(4 + 
      \frac{0.29}{\{1 - 8a(\tilde{v}, \psi)^2\}^2} \right)a(\tilde{v}, \psi).
  \end{align}
\end{lemma}

Note the bound decreases approximately linearly if $a(\tilde{v}, \psi) \to 0$. Finding an exact expression for $a(\tilde{v}, \psi)$ in Lemma
\ref{lem:zhang-ish} is complicated in general, especially since $\tilde{v}$
depends on $\mathcal{I}(\psi)$. However, an upper bound can
often be established, leading to a bound on the density
difference in \eqref{eq:zhang-ish} that is uniform in $v$.

\begin{lemma} \label{lem:zhang-ish-uniform}
  If, for a fixed $\psi \in \mathbb{P}$ and $\bar{a}^2 < 1/8$, it holds for every $v \in \mathbb{S}^{r - 1}$ that  $a(v, \psi) \leq \bar{a}$, then $\vert g(t; v,
  \psi) - \phi(t)\vert \leq 0.14\{4 + 0.29(1 - 8 \bar{a}^2)^{-2}\} \bar{a}$ for
  every $t \in \R{}$ and $v \in \mathbb{S}^{r - 1}$.
\end{lemma}

The lemma allows $\bar{a}$ to depend on $\psi$, but the uniform results we state later result from finding $\bar{a}$ that work for all $\psi$ in some set of interest. Since
$W^S(\psi)$ would be multivariate standard normal under $\psi$ if and only if
$g(t; v, \psi) = \phi(t)$ for every $t\in \R{}$ and $v\in \mathbb{S}^{r - 1}$,
an interpretation of Lemma \ref{lem:zhang-ish-uniform} is that $W^S(\psi)$ is
close to multivariate normal under $\psi$ if $\bar{a}$ is small. By the
arguments following \eqref{eq:model_reml}, Lemma \ref{lem:zhang-ish-uniform} can
be applied to the restricted likelihood, which we give examples of shortly.

To state asymptotic results, we suppose that for each $n \in \{1, 2, \dots\}$, a version of \eqref{eq:marginal} holds. Specifically, $X \in \R{n\times p_n}$, $Z \in \R{n\times q_n}$, $\beta = \beta_n \in \R{p_n}$, $\psi = \psi_n \in \mathbb{P} = \mathbb{P}_n = \{\psi \in \R{r_n}: \Psi_n(\psi)\geq 0, \psi_r > 0\}$, and $\Psi(\cdot) = \Psi_n(\cdot)$ can depend on $n$ in any way consistent with \eqref{eq:marginal}. We omit indices on $Y$, $X$, and $Z$ for simplicity. Let $a_n(\cdot, \cdot)$ and $W^S_n(\cdot)$ be defined as $a(\cdot, \cdot)$ and $W^S(\cdot)$, respectively, with the dependence on $n$ made explicit. 

\begin{theorem} \label{thm:asy_norm}
  Assume that, for $n \in \{1, 2, \dots\}$, \eqref{eq:marginal} holds with known $\beta = 0$ and $\psi = \psi_n \in \mathbb{P}_n \subseteq \R{r_n\times r_n}$. If there exists $\bar{a}_n$ with $\sup_{v \in \mathbb{S}^{r_n - 1}} a_n(v, \psi_n) \leq \bar{a}_n$ for every $n$ and $\lim_{n\to \infty}\bar{a}_n = 0$, then, for any $v_n \in \mathbb{S}^{r_n - 1}$, $v_n^\tsp W_n^S(\psi_n) \to \rN(0, 1)$. Moreover, if in addition $X$ has full column rank for every $n$, and $\beta = \beta_n \in \R{p_n\times p_n}$ is unknown, it holds for any $u_n \in \mathbb{S}^{p_n + r_n
  - 1}$ that $u_n^\tsp W_n^S(\theta_n)\to \rN(0, 1)$.
\end{theorem}
Like Lemma \ref{lem:zhang-ish-uniform}, Theorem \ref{thm:asy_norm} holds if $\bar{a}_n$ depends on $\psi_n$, and the result with known $\beta = 0$ can be applied to the restricted likelihood. The asymptotic normality can hold even if $r_n
\to \infty$, indicating a normality approximation can be useful
even if there are many random effect parameters. What sample size is sufficient in practice will depend on
the setting, including in particular the dependence between observations. The part of the theorem about asymptotic normality when $\beta_n$ is unknown is not immediate from Lemma
\ref{lem:zhang-ish-uniform}. The result is suggested, but not
implied, by the fact that the score for $\beta_n$ is multivariate normal for every $n$. 

We next apply the results of this section in two common settings.

\subsection{Independent clusters} \label{sec:indep_cluster}
Many existing results for mixed models assume a large number of independent
clusters. We consider such a setting here, first stating finite-sample results and then
asymptotic ones. The first result is a separable bound on $a(v,
\psi)$, which does not assume a particular version of
\eqref{eq:marginal} but is often useful when there are many independent
observations or $\Sigma(\psi)$ is well-conditioned. 
 
\begin{lemma} \label{lem:well_cond_bound}
  For any $v \in \R{r}$ and $\psi \in \mathbb{P}$ such that $\Vert \Sigma(v)\Vert_F \neq 0$,
  \begin{align*}
    a(v, \psi) \leq\Vert \Sigma(\psi)^{-1}\Vert \Vert \Sigma(\psi)\Vert
      \frac{\Vert \Sigma(v)\Vert}{\Vert \Sigma(v) \Vert_F}.
  \end{align*}
\end{lemma}
Lemma \ref{lem:well_cond_bound} lets one work separately with $v$ and $\psi$ to
find bounds on $a(v, \psi)$.  For example, when there are $m$
independent clusters of bounded size, then the first product in Lemma
\ref{lem:well_cond_bound} is often bounded uniformly over certain sets of $\psi$, while the ratio is of order $m^{-1/2}$, uniformly in $v \in \mathbb{S}^{r - 1}$.

To be more specific, suppose, independently
for $i \in \{1, \dots, m\}$, for non-stochastic $X_i \in \R{n_i \times p}$ and $Z_i \in \R{n_i \times q_1}$,
\begin{equation} \label{eq:model_clust}
  Y_{i} = X_i\beta + Z_i U_i + E_i,
\end{equation} 
where $E_i \sim \rN(0, \psi_r I_{n_i})$, $U_i \sim \rN(0,
\Psi_1)$, and $\Psi_1 \in \R{q_1\times q_1}$. Now \eqref{eq:marginal} holds for $Y = [Y_1^\tsp, \dots, Y_m^\tsp]^\tsp$ with $X = [X_1^\tsp, \dots, X_m^\tsp]^\tsp \in \R{n\times p}$, and $Z = \bdiag(Z_1, \dots, Z_m) \in \R{n\times q}$, where $\bdiag$ evaluates to a block-diagonal matrix with the arguments as
diagonal blocks. Note $n = \sum_{i = 1}^m n_i =
m\bar{n}$ and $q = mq_1$. We get $\Psi = I_m \otimes \Psi_1$, where $\otimes$ is the Kronecker product, and
\begin{equation} \label{eq:sigma_clust}
  \Sigma = \bdiag(Z_1 \Psi_1 Z_1^\tsp, \dots, Z_m \Psi_1 Z_m^\tsp) + \psi_r I_n.
\end{equation}
Let $\psi_{-r} = \vech(\Psi_1) \in \R{q_1(q_1 + 1)/2}$, the
half-vectorization of $\Psi_1$ obtained by stacking its lower triangular part. Thus, $r = q_1 (q_1 + 1)/2 + 1$ and $\mathbb{P} = \{\psi \in \R{r}: \Psi_1(\psi) \geq 0, \psi_r > 0\}$, where $\Psi_1(\cdot)$ maps $\psi \in \R{r}$ to the $q_1\times q_1$ symmetric matrix whose half-vectorization is $\psi_{-r}$.

To state a result also for the restricted likelihood, define $\tilde{a}(\cdot, \cdot)$ as $a(\cdot, \cdot)$, replacing the $A_j$ with the $\tilde{A}_j$ defined following \eqref{eq:model_reml}. Let also $c_1, \dots, c_4$ denote arbitrary constants that do not
depend on any model parameters or quantities.
\begin{theorem} \label{thm:indep}
  Assume \eqref{eq:model_clust} with (i) $\bar{n}/q_1
  \geq c_1 > 1$ and (ii), for every $i \in \{1, \dots, m\}$,
  \[
    c_2^{-1} \leq \emin(Z_i^\tsp Z_i) \leq \emax(Z_i^\tsp Z_i) \leq c_2,~~c_2 \in (1, \infty).
  \]
 Then for a $c_3 < \infty$,
 \begin{equation} \label{eq:indep_bound}
  a(v, \psi) \leq c_3 m^{-1/2}(1 + \Vert \psi_{-r}\Vert/\psi_r).
 \end{equation}
 Moreover, if in addition (iii) $p/m \leq c_4$ for a small enough $c_4 > 0$, then \eqref{eq:indep_bound} holds with $\tilde{a}(\cdot, \cdot)$ in place of $a(\cdot, \cdot)$, with a different $c_3$.
\end{theorem}

Explicit
expressions for $c_3$ and $c_4$, as functions of $c_1$ and $c_2$, are available
in the proof, but they are somewhat unwieldy. In the special case that $Z_i = Z_1$ for all $i$, we can take $c_3 =  2^{1/2}c_2$.

The $m^{-1/2}$ is intuitive in this setting as $m$ is the number of
independent clusters while observations can be arbitrarily dependent within
clusters. Better bounds are possible by restricting within-cluster
dependence, but the details will depend on the dependence structure.

That $\psi$ appears only through $\Vert \psi_{-r}\Vert/\psi_r$ indicates, by Lemma \ref{lem:zhang-ish-uniform}, that a normality approximation for $W^S(\psi)$ is useful as long as random effect variances are small in comparison to $m^{1/2}$ times the error variance. In particular, it is not a
problem if $\Vert \psi_{-r}\Vert$ is small, including as an extreme special case
$\psi_{-r} = 0$, which corresponds to testing the existence of any random
effects. More generally, testing random effect variances equal to zero or
creating confidence regions for small variances is unproblematic using the
proposed test-statistic. Similar arguments apply to correlations near unity.

It is common to reparameterize the model as $\Sigma
= \sigma^2\{Z \Psi(\tau)Z^\tsp + I_{n}\}$, where $\sigma^2 = \psi_r$ and $\tau =
\psi_{-r}/\psi_r$ \citep[see for
example][]{crainiceanu2004likelihood,bates2005fitting}. For that
parameterization, Theorem \ref{thm:indep} indicates reliable inference is
possible even if $\sigma^2$ is close to zero, as long as $\Vert \tau\Vert$ is
bounded.

The condition that the eigenvalues of $Z_i^\tsp Z_i$ be bounded from above can
be ensured in practice without changing the model, for example by replacing $Z_i$ by $Z_i / \Vert
Z\Vert$. The lower bound on the eigenvalues rules out, for example, the
possibility that some clusters are affected only by a strict subset of the
random effects, as that would mean some $Z_i$ have columns of zeros. This could
be allowed by instead restricting the proportion of clusters effected by each
random effect to be bounded away from zero. However, doing so in full generality
would substantially complicate notation, and we do not believe it would
lead to any fundamental insights.

Theorem \ref{thm:indep} has the following asymptotic result as a corollary. We
omit the proof because it is almost immediate from Theorems \ref{thm:asy_norm}
and \ref{thm:indep}. Recall the notation of Theorem \ref{thm:asy_norm}, but let
us here index by $m$ instead of $n$. Note $m\to \infty$ implies $n\to \infty$
since $n = \sum_{i = 1}^m n_i \geq q_1 m \geq m$.

\begin{corollary} \label{corol:indep_asy}
  Assume that, for $m \in \{1, 2, \dots\}$, \eqref{eq:model_clust} holds with $\beta = \beta_m \in \R{p_m}$ and $\psi = \psi_m \in \mathbb{P}_m \subseteq \R{r_m\times r_m}$. If (i) and (ii) of Theorem \ref{thm:indep} hold for every $m$, with $c_1, c_2$ not depending on $m$; (iv)  $X \in
  \R{n\times p}$ has full column rank for every $m$; and (v) $\Vert \psi_{m}\Vert/\psi_{mr} = o(m^{1/2})$ as
  $m\to \infty$; then it holds for any $u_m \in \mathbb{S}^{r_m + p_m - 1}$ that
  $u_m^\tsp W^S_m(\theta_m) \to \rN(0, 1)$. Moreover, if in addition (iii') $p_m/m \to 0$, then for any $v_m \in
  \mathbb{S}^{r_m - 1}$, $v_m^\tsp \tilde{W}^S_m(\psi_m) \to \rN(0,
  1)$.
\end{corollary}

Quantities in the definition of \eqref{eq:model_clust} that do not appear explicitly in Corollary \ref{corol:indep_asy} can depend on $m$ in any way consistent with \eqref{eq:model_clust} and the assumptions of the corollary. For example, $q_1$ can grow with $m$ as long as (i) holds for every $m$.

In settings where the $n_i$ grow, which is not necessary for Corollary \ref{corol:indep_asy}, it may be more natural to state
bounds for $Z_i^\tsp Z_i / n_i$, but as noted above, $Z_i$ can be standardized
to ensure the upper bound holds. Then, if $n_i \to \infty$, the bounds
effectively say the standardized $Z_i^\tsp Z_i$ should be well-conditioned, for
example by tending some positive definite limit. Such a condition often holds if
the $Z_i$ are drawn randomly, but it rules out, for example, within-cluster
crossed random effects.

If $r_m = r$ is fixed, Corollary \ref{corol:indep_asy}, through Lemma \ref{lem:same_asy_dist}, implies asymptotically correct uniform
coverage probability on compact sets; the key observation is that the conditions of Corollary \ref{corol:indep_asy} are compatible with any sequence $(\psi_m)$ convergent in $\mathbb{P}$.
In this setting, $r$ being fixed implies $q_1$ is also fixed.

We state a result for a restricted score confidence region, which does not require $p_m$ to be fixed, but a
similar result holds for the score confidence region for $\theta$ with unknown $\beta$ if $p_m = p$ is also
fixed. Recall that $c_{r, 1 - \alpha}$ denotes the $(1 - \alpha)$th quantile of $\chi^2_r$, and define $\tilde{\mathbb{C}}_n^S(\alpha)$ by \eqref{eq:region} with $T(\psi) = \tilde{T}_
n^S(\psi)$ and $q_{1 - \alpha}(\psi) = c_{r, 1 - \alpha}$.

\begin{corollary}\label{corol:indep_cover}
  Assume conditions (i), (ii), (iii') and (iv) of Corollary \ref{corol:indep_asy}. Then if $r_m = r$ is fixed as $m\to \infty$, $\tilde{\mathbb{C}}_n^S(\alpha)$
  satisfies \eqref{eq:asy_uni_cover} for any $\alpha \in (0, 1)$
\end{corollary}

\subsection{Crossed random effects} \label{sec:crossed}

With crossed random effects it is impossible to split the data into several
independent vectors, and Lemma \ref{lem:well_cond_bound} is typically not
useful because $\Sigma(\psi)$ is not well-conditioned. Specifically, even for
fixed $\psi$, some eigenvalues usually grow without bound as $n\to \infty$ while
others are bounded. Thus, more intricate arguments are needed to get relevant
asymptotic theory.

To understand the setting, suppose there is one mean parameter $\beta \in \R{}$, and only two random effects, which
is common \citep{jiang2013subset,ekvall2020consistent,jiang2024precise,lyu2024increasing}. Specifically,
for $i \in \{1, \dots, n_1\}$ and $j \in \{1, \dots, n_2\}$,
\begin{equation*} 
  Y_{ij} = \beta + U_{1i} + U_{2j} + E_{ij},
\end{equation*}
where, independently, $U_{i1}\sim \rN(0, \psi_1)$ and $U_{2j} \sim \rN(0,
\psi_2)$, and $E_{ij} \sim \rN(0, \psi_3)$. Let $Y = [Y_{11}, Y_{12}, \dots, Y_{n_1n_2}]^\tsp$, $Z^{(1)} = I_{n_1}
\otimes 1_{n_2} \in \R{n \times n_1}$, $Z^{(2)} = 1_{n_1} \otimes I_{n_2} \in
\R{n \times n_2}$, and $Z = [Z^{(1)}, Z^{(2)}]$, where $1_{n_j} \in \R{n_j}$ is
a vector of ones. Thus, $Z^{(1)}$ and $Z^{(2)}$ correspond to the first and
second random effect, respectively. Now, with $U = [U_{11}, U_{12},
\dots, U_{n_1 n_2}]^\tsp \in \R{n_1 + n_2}$ and $\Psi = \bdiag(\psi_1 I_{n_1},
\psi_2 I_{n_2})$,
\[
  \Sigma = Z\Psi Z^\tsp + \psi_3 I_n = \psi_1 Z^{(1)}Z^{(1)\tsp} 
    + \psi_2 Z^{(2)}Z^{(2)\tsp} + \psi_3 I_n,
\]
where $Z^{(1)}Z^{(1)\tsp} = I_{n_1}\otimes 1_{n_2}1_{n_2}^\tsp$ and
$Z^{(2)}Z^{(2)^\tsp} = 1_{n_1}1_{n_1}^\tsp \otimes I_{n_2}$. One reason crossed random effects are complicated is that $\Sigma$
is not block-diagonal, and cannot be made so by reordering observations. However, a key observation for the subsequent analysis
is that $\Sigma$ is a linear combination of three projection matrices, neither of which depends on $\psi$. Thus, the eigenvalues of $\Sigma$ depend on $\psi$, but the eigenvectors do not. Specifically, since $P_j = 1_{n_j}1_{n_j}^\tsp / n_j$, $j \in \{1, 2\}$, is a projection matrix, so are $\mathcal{P}_1 = Z^{(1)}Z^{(1)\tsp}/n_2 = I_{n_1} \otimes P_{2}$ and $\mathcal{P}_2 = Z^{(2)}Z^{(2)\tsp}/n_2 = P_1 \otimes I_{n_2}$. With these definitions, $\Sigma = \psi_1 n_2 \mathcal{P}_1 + \psi_2 n_1 \mathcal{P}_2 + \psi_3 I_n$.

Similar arguments apply when there are $r - 1 \geq 2$ crossed random effects. Let $Z^{(j)} = 1_{n_1}
\otimes \cdots \otimes 1_{n_{j - 1}} \otimes I_{n_j} \otimes 1_{n_{j + 1}}
\otimes \cdots \otimes 1_{n_{r - 1}}$, $j \in \{1, \dots, r - 1\}$, and $Z =
[Z^{(1)}, \dots, Z^{(r - 1)}]$. Define
$\mathcal{P}_j = P_1 \otimes \cdots \otimes P_{j - 1} \otimes I_{n_j} \otimes
P_{j + 1} \otimes \cdots \otimes P_{r - 1}$, $j \in \{1, \dots, r - 1\}$. Then
with $\Psi = \bdiag(\psi_1 I_{n_1}, \dots, \psi_{r - 1}I_{r - 1})$,
\begin{align} \label{eq:crossed_sigma}
 Y \sim \rN(1_n \beta, \Sigma), ~~ \Sigma = \sum_{j = 1}^{r - 1} \psi_j n_{(j)}\mathcal {P}_j + \psi_r I_{n},
\end{align}
where $n_{(j)} = \prod_{i\neq j} n_i$ and  $n = \prod_{i = 1}^{r - 1} n_i$. To state the main result for crossed random
effects, let also $n_{\min} = \min_{1\leq j\leq r - 1}n_j$ and $\tilde{n} = n - 1 -
\sum_{j = 1}^{r - 1}(n_j - 1)$.

\begin{theorem} \label{thm:crossed}
 Assume \eqref{eq:crossed_sigma}. If $r \geq 3$ and $n_{\min} \geq 2$, then
 \[
  a(v, \psi)^2 \leq \sum_{j = 1}^{r - 1}\frac{1}{n_j - 1} + \frac{(r - 2)^2}{\tilde{n}}
   \leq \frac{r - 1}{n_{\min} - 1} + \frac{(r - 2)^2}{\tilde{n}}.
 \]
 Moreover, if also $n_{\min} \geq 3$, then
 \[
    \tilde{a}(v, \psi)^2 \leq \max\left\{(n_{\min} - 2)^{-1},  (\tilde{n} - 1)^{-1}\right\}.
 \]
\end{theorem}

The first part of the proof of Theorem \ref{thm:crossed} consists of finding a
spectral decomposition of $\Sigma$ (Supplementary Material). We provide a direct
argument but note that the results of \citet{henderson1981deriving}, which were
used by \citet{lyu2024increasing}, could also be used for that particular part of the
proof. We prefer the direct argument because it is short and sets up the rest of
the proof.

The bounds in Theorem \ref{thm:crossed} are remarkably simple. That they do not
depend on $\psi$ indicates reliable inference is possible uniformly over the
parameter set. It is also notable that the bound for the restricted likelihood
can be better, even if $\beta = 0$ is known. The technical reasons for this can
be seen in the proof. The intuition is that centering can facilitate variance
estimation even if the population mean is known to be zero. We only considered a
simple mean structure here because our focus is on $\psi$. We expect similar
results are possible for more general means, for example by using the idea of
\citet{lyu2024increasing} to decompose the predictors into parts varying only along
certain dimensions.

What is required in general is that $n_{\min}$ is large in comparison to
$r$. When this holds, $\tilde{n}$ will also be large in comparison to $r^2$ under
mild conditions. Thus, the result indicates the complicated dependence induced
by crossed random effects decreases the effective number of observations from
$n$ to $n_{\min}$. For example, increasing $r$ while keeping $n_{\min}$ fixed
will increase $n$ but may worsen the bound. 

Corollaries with asymptotic results similar to those in Section
\ref{sec:indep_cluster} are straightforward from Theorem \ref{thm:crossed}. For
brevity, we only state one for the restricted score confidence region, $\tilde{\mathbb{C}}_n^S(\alpha)$. The result is remarkable because the uniformity is over the entire parameter
set; this is an effect of the bounds in Theorem \ref{thm:crossed} not depending on $\psi$. Because $n = \prod_{j = 1}^{r - 1} n_j$, which can only take certain integer values, we index by $k \in \{1, 2, \dots\}$. As in previous asymptotic results, model quantities not appearing explicitly in the statement of the corollary can depend on $k$ in any way consistent with the assumptions of the corollary.

\begin{corollary} \label{corol:cross}
  Assume, for each $k \in \{1, 2, \dots\}$, \eqref{eq:crossed_sigma} holds with an $n_{\min} = n^k_{min}$ such that $n^{k}_{\min} \to \infty$ as $k\to \infty$, and fixed $r \geq 3$. Then
  \[
    \sup_{\psi \in [0, \infty)^{r - 1}\times (0, \infty)}\left \vert \pr_{\psi}\{\psi \in 
    \tilde{\mathbb{C}}_n^S(\alpha)\} - (1 - \alpha)\right\vert \to 0.
  \]
\end{corollary}

  \section{Numerical experiments} \label{sec:sims}

  We investigate fine-sample coverage probabilities for confidence regions for
  $\psi$ in three different settings. In the first, we generate data
  from a version of the model with independent clusters and correlated random effects discussed in Section
  \ref{sec:indep_cluster}. For $i \in \{ 1, \dots, m\}$ and $j \in \{1, 2, 3\}$,
  \[
    Y_{ij} = X_{ij}^\tsp \beta + Z_{ij}^\tsp U_i + E_{ij},
  \]
  where, independently for all $i$ and $j$, $E_{ij} \sim \rN(0, \psi_4)$ and $U_i
  \sim \rN(0, \Psi_1)$ with
  \[
    \Psi_1 = \begin{bmatrix} \psi_1  & \psi_2 \\ \psi_2 & \psi_3\end{bmatrix}.
  \]
  We fix $\psi_1 = \psi_3 = \psi_4 = 1$ in the simulations and consider different values of $\psi_2$, which is then the correlation between the random effects. We consider the restricted score test-statistic, $\tilde{T}^S(\psi)$, and the corresponding confidence region, $\tilde{\mathbb{C}}^S(\alpha)$. For comparison, we include a confidence region based on the statistic $T^P(\psi) = S(\psi; \tilde{\beta})^\tsp
  \mathcal{I}(\psi)^{-1}S(\psi; \tilde{\beta})$, where $\tilde{\beta} = \{X^\tsp
  \Sigma(\psi)^{-1}X\}^{-1}X^\tsp \Sigma(\psi)^{-1} Y$, which is in effect using the profile score (Supplementary Material). We also include confidence regions based on the
  restricted likelihood ratio and a Wald statistic using
  restricted maximum likelihood estimates standardized by the Fisher information
  at those estimates. We focus on the restricted likelihood because it tends to give more accurate coverage probabilities for all considered methods.
  
  Before starting the simulation, we drew the elements of $\beta$
  from a standard normal distribution, and created the $X_{ij} = [X_{ij1},
  \dots, X_{ijp}]^\tsp \in \R{p}$ by setting $X_{ij1} = 1$ and drawing the
  remaining $X_{ijk}$ independently from a uniform distribution on $(-1, 1)$. We
  set $Z_{ij} = [X_{ij1}, X_{ij2}]^\tsp \in \R{2}$, so there is a random
  intercept and a random slope.

\begin{figure}
    \includegraphics[width=\linewidth]{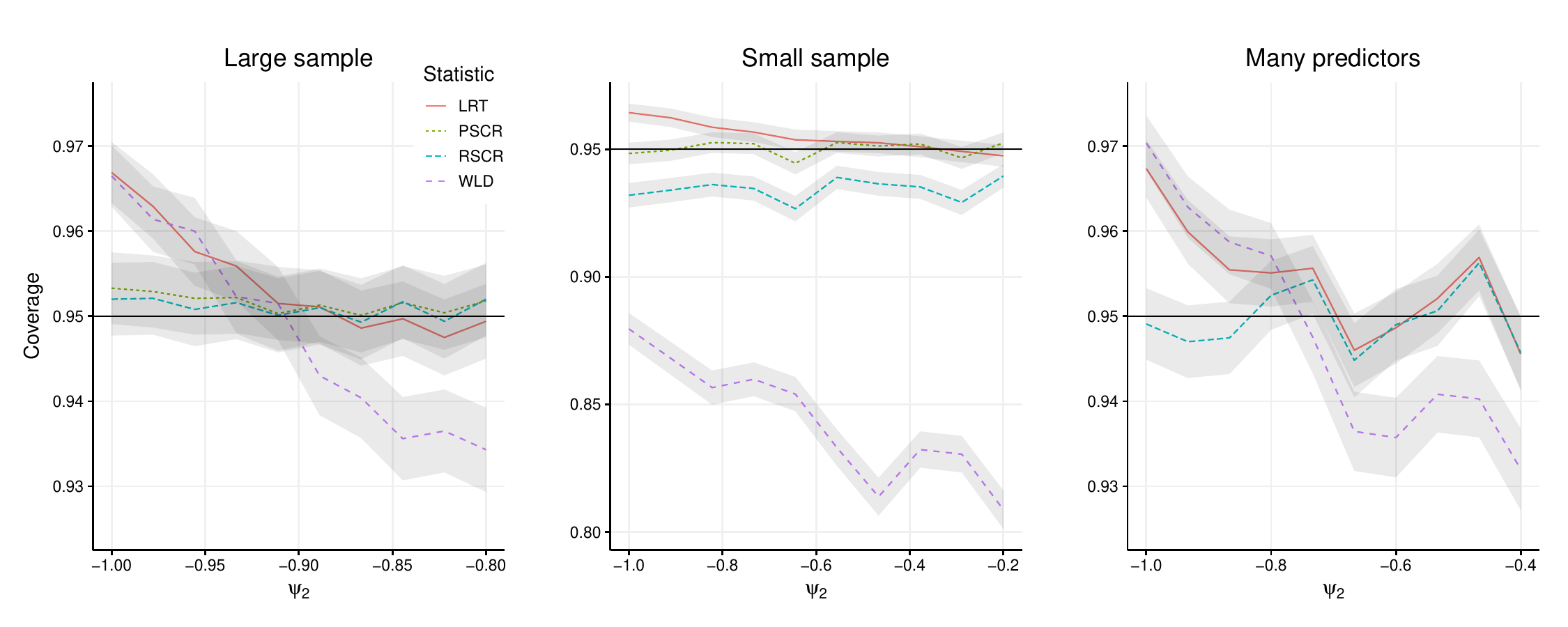}
    \caption{Coverage probabilities with independent clusters and correlated
    random effects, for different values of the random effect correlation $\psi_2 (\psi_1 \psi_3)^{-1/2} = \psi_2$ (horizontal axes).
    The settings are $(m, n_i, p)$ equal to $(500, 3, 2)$ (left plot), $(20, 3,
    2)$ (middle plot), and $(200, 3, 100)$ (right plot).} \label{fig:corr}
\end{figure}

  Figure \ref{fig:corr} shows coverage probabilities for different values of the
  random effect correlation, $\psi_2$ (horizontal axes), for three different
  choices of $m$ and $p$. In the left plot, $m = 500$, so $n =
  1500$, and $p = 2$. As suggested by our theory, the restricted score regions (RSCR)
  have approximately correct coverage probability regardless of how close
  $\psi_1$ is to the boundary point $-1$. The regions based on the profile score
  (PSCR) also have good coverage properties. By contrast, the coverage
  probabilities for the Wald (WLD) and likelihood ratio (LRT) confidence regions depend
  substantially on the true $\psi_2$.
  
  When the sample size is small (middle plot, Figure \ref{fig:corr}), the Wald
  regions are particularly unreliable while the other three are relatively
  similar. The profile score region's coverage is near-nominal for all values of $\psi_2$ while the restricted score region's is below nominal, around 0.92--0.94 for all values of $\psi_2$. The likelihood ratio region also has coverage close to nominal in this setting, with some distortion near the boundary. We caution that this should not be interpreted as the likelihood ratio
  performing well in small samples in general. We do not have a formal result
  for the likelihood ratio region in small samples, but simulations show its coverage
  depends on the setting. What is true more generally is that its coverage near
  the boundary tends to move away from the nominal level as the sample size increases.
  
  Intuitively, there is little difference between restricted and profile score regions when $m$ is large
  in comparison to $p$. Conversely, when $p = 100$ is relatively large in
  comparison to $m = 200$ (right plot, Figure \ref{fig:corr}), the profile
  likelihood has coverage probabilities so far from nominal, around 0.4, that
  we omit it from the plot. In this setting, too, the restricted score gives near-nominal coverage, but the restricted likelihood ratio or Wald statistics do not in general.


\begin{figure}
    \includegraphics[width=\linewidth]{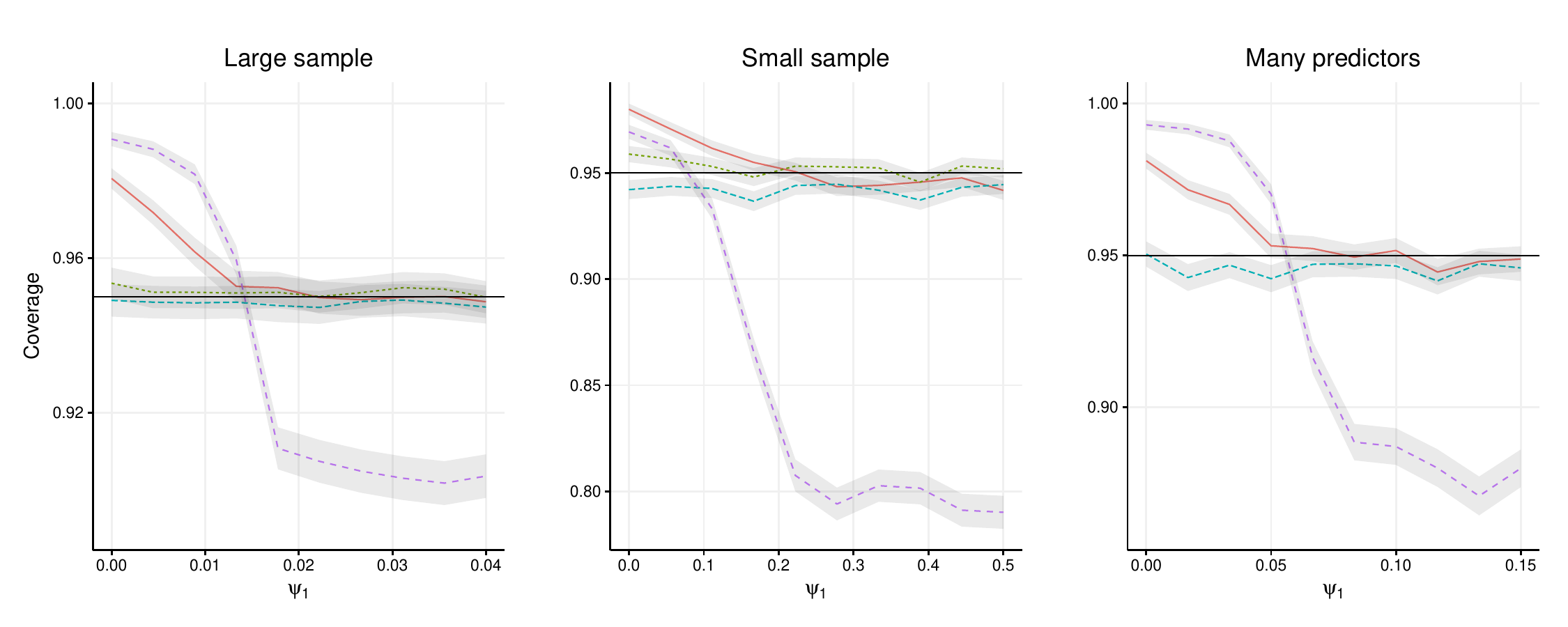}
    \caption{Coverage probabilities with crossed, independent random effects, for
    different random effect variances $\psi_1 = \psi_2$ (horizontal axes). The settings are $(n_1,
    n_2, p)$ equal to $(40, 40, 2)$ (left plot), $(10, 10, 2)$ (middle plot), and
    $(20, 20, 80)$ (right plot).} \label{fig:cross}
\end{figure}

  The second setting is a version of the crossed random effects model discussed in
  Section \ref{sec:crossed}. Specifically, for $i \in \{1, \dots, n_1\}$ and $j
  \in \{1, \dots, n_2\}$,
  \[
      Y_{ij} =  X_{ij}^\tsp \beta + U_{1i} + U_{2j} + E_{ij},
  \]
  with $X_{ij}$, $\beta$, and $E_{ij}$ as in the first setting. The $U_{1i}$ and $U_{2j}$ are independent with variances $\psi_1$ and $\psi_2$, respectively. Data were generated with $\psi_1 = \psi_2$ indicated on the horizontal axis in Figure \ref{fig:cross}. We consider $(n_1, n_2, p)$ equal to $(40, 40, 2)$, so $n = 1600$ is large relative to $p$; $(10, 10, 2)$, so both $n = 100$ and $p$ are relatively small; and $(20, 20, 80)$, so $n = 400$ and $p$ are both moderately large. Recall, because the random effects are crossed, neither of $n$, $n_1$, and $n_2$ is a number of independent
  observations. Nevertheless, Theorem \ref{thm:crossed} suggests $n_{\min} = \min(n_1, n_2)$ controls the quality of the distribution approximations.

  When $n_{\min}$ is relatively large (left plot, Figure \ref{fig:cross}), both score-based confidence region have approximately correct coverage for every value of $\psi_1 = \psi_2$ while the likelihood ratio and Wald regions do not. The Wald regions are particularly unreliable; they are conservative near the boundary and invalid further from the boundary. The difference between the restricted and profile score-based regions is small since $p$ is small.

  When $n_{\min}$ is small (middle plot, Figure \ref{fig:cross}), the profile score region is conservative near the boundary. The restricted score region has coverage slightly below nominal. Nevertheless, both are in general closer to nominal than the likelihood ratio and Wald regions. The likelihood ratio is the most conservative near the boundary, while the Wald region is again conservative near the boundary and invalid further from it.

  When the number of predictors is relatively large (right plot, Figure \ref{fig:cross}), the restricted score-based region has coverage slightly below but near nominal for all values of $\psi_1 = \psi_2$. The likelihood ratio region is conservative near the boundary, but the Wald region is even more so. As in the other settings, the Wald region is invalid further from the boundary.

  The Supplementary Material includes additional simulation settings, and the results are qualitatively similar to those presented here.
  
  \section{Final remarks} \label{sec:final}
  
  Our results show the score standardized by expected information is
  often approximately multivariate normally distributed, even with complicated
  dependence and a large number of parameters. This leads to reliable inference.
  Such results are impossible for Wald and likelihood ratio statistics when
  maximum likelihood estimators can be on the boundary with appreciable
  probability. Nevertheless, asymptotically uniformly correct inference with those statistics
  may be possible even in the settings we consider, but how to achieve that is
  currently unclear and an avenue for future research.
  
  We expect our results can be extended to inference on other types of
  covariance matrices in linear models or, more generally, inference using
  multivariate normal likelihood as a pseudo-likelihood for possibly non-normal
  data. Non-normally distributed random effects or errors could be allowed by
  developing a distributional approximation similar to Lemma
  \ref{lem:zhang-ish-uniform} for non-normal quadratic forms. One would also
  need to address the fact that the covariance matrix of the score depends on
  the third and fourth moment of the responses, possibly by using resampling or
  Monte Carlo-based methods. Similar asymptotic results could likely be obtained
  without assuming normality, either using Lindeberg's conditions directly on
  the quadratic forms or by adapting more specialized results \citep[see][for
  example]{jiang1996reml}. For settings with crossed random effects, it may be
  possible to adapt some of the asymptotic results for the score function in
  \citet{lyu2024increasing} and \citet{jiang2024precise} to boundary settings. However,
  substantial work may be needed to make them appropriately uniform.
  
  The proposed confidence regions can be computed efficiently using the
  expressions for the score and information matrix in Section \ref{sec:model}. Some additional remarks on implementation are in the Supplementary Material. Code for reproducing the simulation results is at \url{https://github.com/koekvall/uniform_lmm_suppl/}. 
  
  It is in principle possible to extend some of our results to generalized linear mixed
  models. However, both theory and computation will be substantially more
  difficult, especially for the case with crossed random effects. Indeed, in
  that setting even pointwise asymptotic theory with interior parameters is
  challenging \citep{jiang2025asymptotic}.
  
  Finally, while we have not focused on the testing on boundary points because
  that is a well-studied problem, we note tests with asymptotically correct size
  are immediate from our more general results.

\section*{Acknowledgement}
The authors are grateful to two referees and an Associate Editor for comments that improved the manuscript substantially. The authors also thank James Hodges and Aaron Molstad for helpful discussions, and Yiqao Zhang and Matias Shedden for assisting with programming.

\section*{Supplementary material}
\label{SM}
Proofs of lemmas and theorems are in the Supplementary Material along with additional numerical results.  

\bibliographystyle{apalike}
\bibliography{uniform_lmm}

\end{document}